# Approximate Solution of Kuramoto–Sivashinsky Equation Using Reduced Differential Transform Method

Omer ACAN[1,a] and Yıldıray KESKİN[1,b]

*Department of Mathematics, Science Faculty, Selcuk University, Konya 42003, Turkey.*

[a]acan_omer@ selcuk.edu.tr
[b]yildiraykeskin@yahoo.com

**Abstract.** In this study, approximate solution of Kuramoto–Sivashinsky Equation, by the reduced differential transform method, are presented. We apply this method to an example. Thus, we have obtained numerical solution Kuramoto–Sivashinsky equation. Comparisons are made between the exact solution and the reduced differential transform method. The results show that this method is very effective and simple.

## 1. INTRODUCTION

Partial differential equations (PDEs) have numerous essential applications in various fields of science and engineering such as fluid mechanic, thermodynamic, heat transfer, physics [10].

Most new nonlinear PDEs do not have a precise analytic solution. So, numerical methods have largely been used to handle these equations. It is difficult to handle nonlinear part of these equations. Although most of scientists applied numerical methods to find the solution of these equations, solving such equations analytically is of fundamental importance since the existent numerical methods which approximate the solution of PDEs don't result in such an exact and analytical solution which is obtained by analytical methods.

Many researchers have paid attention by studying to the solutions of nonlinear PDEs by various methods [11-19]. RDTM [1-4], devised by Yıldıray Keskin in 2009, is a numerical method to obtain approximate solutions of various types of nonlinear partial differential equations. It has received much attention since it has applied to solve a wide variety of problems by many authors [20-23].

In this study, RDTM is used to obtain approximate solution of KS equation. The Generalized Kuramoto-Sivashinsky equation is given [5,6] by

$$u_t + \alpha u^\beta u_x + \gamma u^\tau u_{xx} + \lambda u_{xxxx} = 0, \tag{1.1}$$

Where $\alpha, \beta, \gamma, \lambda, \tau \neq 0$ and $\alpha, \beta, \gamma, \lambda, \tau \in R$.

When $\alpha = \beta = 1$ and $\tau = 0$, eq. (1.1) reduces to original KS equation of the form [6-8]

$$u_t + uu_x + \gamma u_{xx} + \lambda u_{xxxx} = 0. \tag{1.2}$$

In the past several decades, many researchers have used various methods to solve KS equation. We apply the RDTM to solve KS equation of the form (1.2).

## 2. ANALYSIS OF THE RDTM

The basic definitions in the RDTM [4] are as follows:

**Definition 2.1.**

If function $u(x,t)$ is analytic and differentiated continuously with respect to time t and space x in the domain of interest, then let

$$U_k(x) = \frac{1}{k!}\left[\frac{\partial^k}{\partial t^k}u(x,t)\right]_{t=0} \qquad (2.1)$$

where the t-dimensional spectrum function $U_k(x)$ is the transformed function. In this paper, the lowercase $u(x,t)$ represent the original function while the uppercase $U_k(x)$ stand for the transformed function.

The differential inverse transform of $U_k(x)$ is defined as follows:

$$u(x,t) = \sum_{k=0}^{\infty} U_k(x) t^k. \qquad (2.2)$$

Then combining equation (2.1) and (2.2) we write

$$u(x,t) = \sum_{k=0}^{\infty} \frac{1}{k!}\left[\frac{\partial^k}{\partial t^k}u(x,t)\right]_{t=0} t^k. \qquad (2.3)$$

From the above definitions, it can be found that the concept of the reduced differential transform is derived from the power series expansion.

For the purpose of illustration of the methodology to the proposed method, we write the KS equation in the standard operator form

$$L_t(u(x,t)) + L_x(u(x,t)) + N(u(x,t)) = 0 \qquad (2.4)$$

with initial condition

$$u(x,0) = f(x) \qquad (2.5)$$

where $L_t(u(x,t)) = u_t$ and $L_x(u(x,t)) = \gamma u_{xx} + \lambda u_{xxxx}$ are linear operators which have partial derivatives, $N(u(x,t)) = u u_x$ is a nonlinear term.

**TABLE 1.** The fundamental operators of RDTM [1-4]

| Functional Form | Transformed Form |
|---|---|
| $u(x,t)$ | $U_k(x) = \frac{1}{k!}\left[\frac{\partial^k}{\partial t^k}u(x,t)\right]_{t=0}$ |
| $w(x,t) = u(x,t) \pm \alpha v(x,t)$ | $W_k(x) = U_k(x) \pm \alpha V_k(x)$ ($\alpha$ is a constant) |
| $w(x,t) = x^m t^p u(x,t)$ | $W_k(x) = x^m U_{(k-p)}(x)$ |
| $w(x,t) = u(x,t) v(x,t)$ | $W_k(x) = \sum_{r=0}^{k} V_r(x) U_{k-r}(x) = \sum_{r=0}^{k} U_r(x) V_{k-r}(x)$ |
| $w(x,t) = \frac{\partial^r}{\partial t^r} u(x,t)$ | $W_k(x) = (k+1)\ldots(k+r) U_{k+r}(x) = \frac{(k+r)!}{k!} U_{k+r}(x)$ |
| $w(x,t) = \frac{\partial}{\partial x^m} u(x,t)$ | $W_k(x) = \frac{\partial}{\partial x^m} U_k(x)$ |

According to the RDTM and Table 1, we can construct the following iteration formula:

$$(k+1)U_{k+1}(x) = -\gamma \frac{\partial^2}{\partial x^2} U_k(x) - \lambda \frac{\partial^4}{\partial x^4} U_k(x) - \sum_{r=0}^{k} U_{k-r}(x) \left(\frac{\partial^2}{\partial x^2} U_r(x)\right). \qquad (2.6)$$

From initial condition (2.2), we write

$$U_0(x) = f(x) \qquad (2.7)$$

substituting (2.7) into (2.6) and by a straight forward iterative calculations, we get the following $U_k(x)$ values. Then the inverse transformation of the set of values $\{U_k(x)\}_{k=0}^{n}$ gives approximation solution as,

$$\tilde{u}_n(x,t) = \sum_{k=0}^{n} U_k(x) t^k \tag{2.8}$$

where *n* is order of approximation solution.

Therefore, the exact solution of problem is given by
$$u(x,t) = \lim_{n \to \infty} \tilde{u}_n(x,t). \tag{2.9}$$

## 3. NUMERICAL APPLICATIONS

In this section, we test the RDTM for KS equation.

**Example 3.1.** We consider the following KS equation [6,9]:
$$u_t + u u_x + u_{xx} + u_{xxxx} = 0 \tag{3.1}$$

subject to initial condition
$$u(x,0) = c + \frac{5}{19}\sqrt{\frac{11}{19}} \left(11 \tanh^3\left(k(x-x_0)\right) - 9 \tanh\left(k(x-x_0)\right)\right). \tag{3.2}$$

In [6,9], the exact solution of (3.1) is given as
$$u(x,t) = c + \frac{5}{19}\sqrt{\frac{11}{19}} \left(11 \tanh^3\left(k(x-ct-x_0)\right) - 9 \tanh\left(k(x-ct-x_0)\right)\right). \tag{3.3}$$

Taking differential transform of (3.1) and the initial condition (3.2) respectively, we obtain
$$(k+1) U_{k+1}(x) = -\gamma \frac{\partial^2}{\partial x^2} U_k(x) - \lambda \frac{\partial^4}{\partial x^4} U_k(x) - \sum_{r=0}^{k} U_{k-r}(x) \left(\frac{\partial^2}{\partial x^2} U_r(x)\right). \tag{3.3}$$

where the t-dimensional spectrum function $U_k(x)$ are the transformed function.

From the initial condition (3.2) we write
$$U_0(x) = c + \frac{5}{19}\sqrt{\frac{11}{19}} \left(11 \tanh^3\left(k(x-x_0)\right) - 9 \tanh\left(k(x-x_0)\right)\right). \tag{3.4}$$

Now, substituting (3.4) into (3.3), we obtain the following $U_k(x)$ values successively

$$U_1(x) = -\frac{60}{361} \sinh\left(k(-x+x_0)\right) \begin{pmatrix} 4\cosh^4\left(k(-x+x_0)\right) - 11\cosh^2\left(k(-x+x_0)\right) + 16k^2 \times \\ \cosh^4\left(k(-x+x_0)\right) - 24k^2 \cosh^2\left(k(-x+x_0)\right) + 330k^2 \end{pmatrix} \frac{k^2\sqrt{209}}{\cosh^7\left(k(-x+x_0)\right)} \tag{3.5}$$

$$U_2(x) = -\frac{60}{361} \sinh\left(k(-x+x_0)\right) \begin{pmatrix} 831600k^4 + 165\cosh^4\left(k(-x+x_0)\right) + 18600k^2 \times \\ \cosh^4\left(k(-x+x_0)\right) - 18480k^2 \cosh^2\left(k(-x+x_0)\right) \\ +8k^2 \cosh^8\left(k(-x+x_0)\right) + 391440k^4 \cosh^4\left(k(-x+x_0)\right) \\ -1112160k^4 \cosh^2\left(k(-x+x_0)\right) + 64k^2 \cosh^8\left(k(-x+x_0)\right) \\ -3776k^2 \cosh^6\left(k(-x+x_0)\right) + 128k^4 \cosh^8\left(k(-x+x_0)\right) \\ -30592k^4 \cosh^6\left(k(-x+x_0)\right) - 112\cosh^6\left(k(-x+x_0)\right) \end{pmatrix} \frac{k^4\sqrt{209}}{\cosh^{11}\left(k(-x+x_0)\right)}$$

(3.6)

From (2.8)
$$\tilde{u}_2(x,t) = \sum_{k=0}^{2} U_k(x) t^k. \tag{3.7}$$

Substituting (3.4), (3.5) and (3.6) in (3.7), we have

$$\tilde{u}_2(x,t) = -\frac{1}{361\cosh^{11}(k(-x+x_0))} \begin{pmatrix} -361c\cosh^{11}(k(-x+x_0)) + 10\sqrt{209}\sinh(k(-x+x_0))\cosh^{10}(k(-x+x_0)) - 55\sqrt{209} \times \\ \sinh(k(-x+x_0))\cosh^8(k(-x+x_0)) + 240\sqrt{209}k^2t\sinh(k(-x+x_0))\cosh^8(k(-x+x_0)) \\ -660\sqrt{209}k^2t\sinh(k(-x+x_0))\cosh^6(k(-x+x_0)) + 960\sqrt{209}k^4t\sinh(k(-x+x_0)) \times \\ \cosh^8(k(-x+x_0)) - 13440\sqrt{209}k^4t\sinh(k(-x+x_0))\cosh^6(k(-x+x_0)) + 19800\sqrt{209} \times \\ k^4t\sinh(k(-x+x_0))\cosh^4(k(-x+x_0)) - 49896000\sqrt{209}k^8t^2\sinh(k(-x+x_0)) - 9900 \times \\ \sqrt{209}k^4t^2\sinh(k(-x+x_0))\cosh^4(k(-x+x_0)) - 1116000\sqrt{209}k^6t^2\sinh(k(-x+x_0)) \times \\ \cosh^4(k(-x+x_0)) + 1108800\sqrt{209}k^6t^2\sinh(k(-x+x_0))\cosh^2(k(-x+x_0)) - 480\sqrt{209} \times \\ k^4t^2\sinh(k(-x+x_0))\cosh^8(k(-x+x_0)) - 23486400\sqrt{209}k^8t^2\sinh(k(-x+x_0)) \times \\ \cosh^4(k(-x+x_0)) + 66729600\sqrt{209}k^8t^2\sinh(k(-x+x_0))\cosh^2(k(-x+x_0)) - 3840 \times \\ \sqrt{209}k^6t^2\sinh(k(-x+x_0))\cosh^8(k(-x+x_0)) + 226560\sqrt{209}k^6t^2\sinh(k(-x+x_0)) \times \\ \cosh^6(k(-x+x_0)) - 7680\sqrt{209}k^8t^2\sinh(k(-x+x_0))\cosh^8(k(-x+x_0)) + 1835520\sqrt{209}k^8 \times \\ t^2\sinh(k(-x+x_0))\cosh^6(k(-x+x_0)) + 6720\sqrt{209}k^4t^2\sinh(k(-x+x_0))\cosh^6(k(-x+x_0)) \end{pmatrix}$$

When analysed the solution of the above an example by RDTM, the following results are obtained:

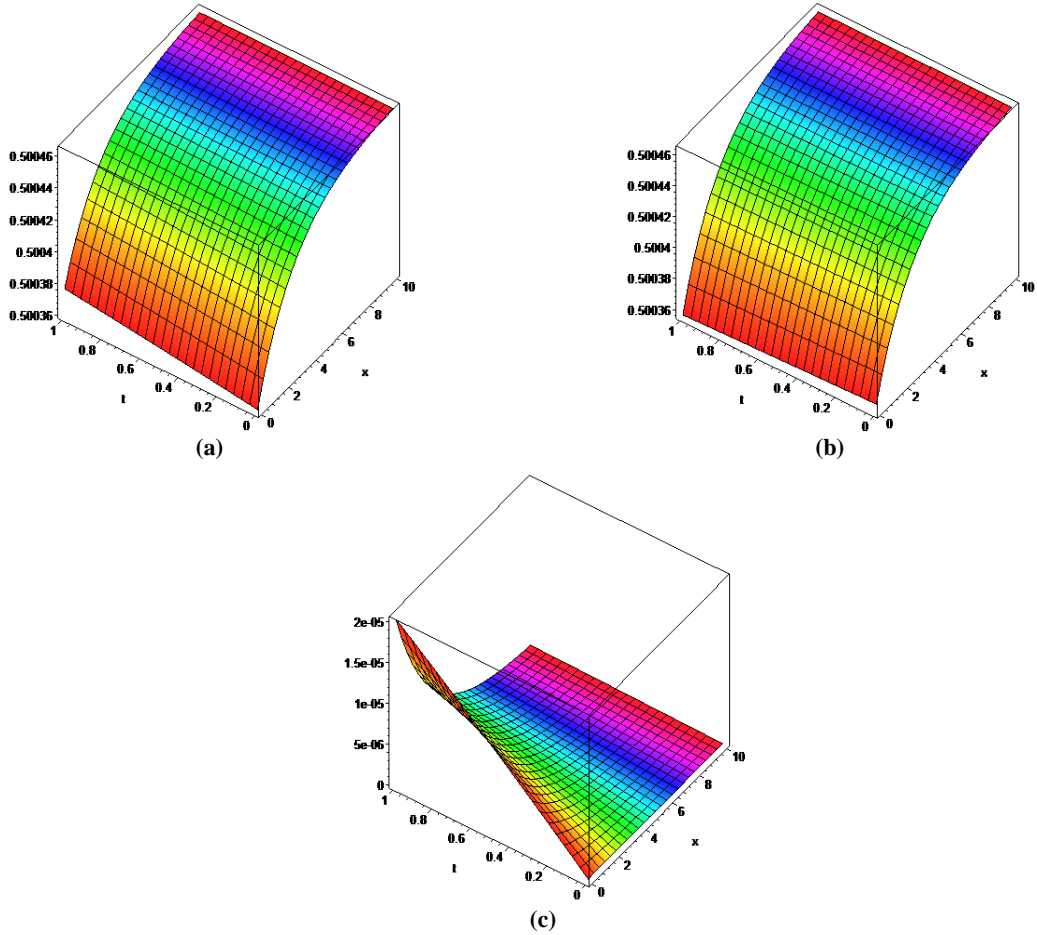

**FIGURE 1.** The surface shows the solution $u(x,t)$ for eq. (3.1) when $c = 0.1$, $k = \left(\sqrt{11/19}\right)/4$, $x_0 = -30$ (a) exact solution (3.3.), (b) 2th order of approximate solution for eq. (3.1) and (c) the absolute error between exact and numerical solution.

**TABLE 1**. The exact and 2th order of approximate solution for (3.1) and absolute error between exact and numerical solutions

| $x$ | $t$ | RDTM Solution | Exact Solution | Abs-Error |
|---|---|---|---|---|
| 0.0 | 0.0 | 0.5003600908 | 0.500360093 | $0.5420022964 \times 10^{-9}$ |
| 0.0 | 0.5 | 0.5003685212 | 0.500358052 | $0.1046850860 \times 10^{-4}$ |
| 0.0 | 1.0 | 0.5003762240 | 0.500355975 | $0.2024957166 \times 10^{-4}$ |
| 0.5 | 0.0 | 0.5003784827 | 0.500378483 | $0.6041276398 \times 10^{-9}$ |
| 0.5 | 0.5 | 0.5003854531 | 0.500376798 | $0.8655395022 \times 10^{-5}$ |
| 0.5 | 1.0 | 0.5003918215 | 0.500375080 | $0.1674210244 \times 10^{-4}$ |
| 1.0 | 0.0 | 0.5003936888 | 0.500393690 | $0.5954475010 \times 10^{-9}$ |
| 1.0 | 0.5 | 0.5003799452 | 0.500392296 | $0.7156111486 \times 10^{-5}$ |
| 1.0 | 1.0 | 0.5004047174 | 0.500390875 | $0.1384174699 \times 10^{-4}$ |

## 4. CONCLUSION

In this study, reduced differential transform method has been successfully applied to KS equation. This example shows that the results of the present method are in excellent agreement with those of absolute errors and the obtained solutions shown in table.